\documentclass[a4paper,11pt]{amsart}
\usepackage{graphicx}
\usepackage{mathptmx}
\usepackage{amsmath}
\usepackage{amssymb}
\usepackage{enumitem}
\usepackage{xcolor}
\usepackage{xparse}
\NewDocumentCommand{\eulerian}{omm}
 {%
  \genfrac<>{0pt}{}{#2}{#3}%
  \IfValueT{#1}{_{\!#1}}%
 }

 \newmuskip\pFqmuskip

\newcommand*\pFq[6][8]{%
  \begingroup 
  \pFqmuskip=#1mu\relax
  \mathchardef\normalcomma=\mathcode`,
  \mathcode`\,=\string"8000
  \begingroup\lccode`\~=`\,
  \lowercase{\endgroup\let~}\pFqcomma
  {}_{#2}F_{#3}{\left(\genfrac..{0pt}{}{#4}{#5}\bigg|#6\right)}%
  \endgroup
}
\newcommand{\pFqcomma}{{\normalcomma}\mskip\pFqmuskip}

\newtheorem{theorem}{Theorem}

\newtheorem{remark}[theorem]{Remark}

\begin{document}

\title[Study on $r$-truncated degenerate Stirling numbers of the second kind]{Study on $r$-truncated degenerate Stirling numbers of the second kind}

\author{Taekyun  Kim*}
\address{Department of Mathematics, Kwangwoon University, Seoul 139-701, Republic of Korea}
\email{tkkim@kw.ac.kr}

\author{DAE SAN KIM*}
\address{Department of Mathematics, Sogang University, Seoul 121-742, Republic of Korea}
\email{dskim@sogang.ac.kr}

\author{Hyekyung Kim*}
\address{Department Of Mathematics Education, Daegu Catholic University, Gyeongsan 38430, Republic of Korea}
\email{hkkim@cu.ac.kr}

\thanks{ * are corresponding authors.}

\maketitle

\begin{abstract}
The degenerate Stirling numbers of the second kind and of the first kind, which are respectively degenerate versions of the  Stirling numbers of the second kind and of the first kind, appear frequently when we study various degenerate versions of some special numbers and polynomials.
The aim of this paper is to consider the $r$-truncated degenerate Stirling numbers of the second kind, which reduce to the degenerate Stirling numbers of the second for $r=1$, and to investigate their explicit expressions, some properties and related identities, in connection with several other degenerate special numbers and polynomials.
\end{abstract}

\section{Introduction}
Carlitz [3,4] initiated to study the degenerate Bernoulli and Euler polynomials and numbers, which are degenerate versions of the Bernoulli and Euler polynomials and numbers.
In recent years, studying degenerate versions of some special numbers and polynomials have regained interests of some mathematicians and yielded quite a few interesting results (see [8-15,17] and the references therein). It is noteworthy that studying degenerate versions is not only limited to polynomials but also extended to transcendental functions, like gamma functions. It is also remarkable that the degenerate umbral calculus is introduced as a degenerate version of the classical umbral calculus. Degenerate versions of special numbers and polynomials have been explored by various methods, including combinatorial methods, generating functions, umbral calculus techniques, $p$-adic analysis, differential equations, special functions, probability theory and analytic number theory.\par

The Stirling number of the second $S_{2}(n,k)$ is the number of ways to partition a set of $n$ objects into $k$ nonempty subsets. The (signed) Stirling number of the first kind $S_{1}(n,k)$ is defined such that the number of permutations of $n$ elements having exactly $k$ cycles is the nonnegative integer $(-1)^{n-k}S_{1}(n,k)=|S_{1}(n,k)|$. The degenerate Stirling numbers of the second kind $S_{2,\lambda}(n,k)$ (see \eqref{5}, \eqref{7}) and of the first kind $S_{1,\lambda}(n,k)$ (see \eqref{4}, \eqref{6-1}) were introduced as degenerate versions of the Stirling numbers of the second and of the first kind, respectively. These degenerate Stirling numbers of both kinds appear frequently when we study degenerate versions of some special numbers and polynomials.\par

Here we consider the $r$-truncated degenerate Stirling numbers of the second which reduce to the degenerate Stirling numbers of the second for $r=1$.
The aim of this paper is by using generating functions to study their explicit expressions, some properties and related identities on the $r$-truncated degenerate Stirling numbers of the second kind, in connection with the higher-order $r$-truncated degenerate Bernoulli polynomials, the degenerate Bernoulli numbers, the degenerate Stirling numbers of the second kind and the higher-order degenerate Bernoulli numbers.\par
The outline of this paper is as follows. In Section 1, we recall the degenerate exponentials and logarithms. We remind of the degenrate Stirling numbers of both kinds. Furthermore, we recall the higher-order degenerate Bernoulli polynomials. Section 2 is the main results of this paper. We consider the $r$-truncated degenerate Stirling numbers of the second kind and find three explicit expressions for those numbers in Theorems 1-3. Then we introduce the higher-order $r$-truncated degenerate Bernoulli polynomials and obtain some results in connection with the $r$-truncated degenerate Stirling numbers of the second kind. In Theorem 4, the degenerate Bernoulli numbers are expressed in terms of the Stirling numbers of the second kind. In Theorems 5 and 8, obtained are some identities involving the $r$-truncated degenerate Stirling numbers of the second kind (for $r$=2 and $r$=3, respectively), the degenerate Bernoulli numbers and the degenerate Stirling numbers of the second. In Theorem 6, we find an identity connecting the 2-truncated degenerate Stirling numbers of the second kind and the degenerate Bernoulli numbers. In Theorem 7, we get an identity relating the 2-truncated degenerate Stirling numbers of the second kind and the higher-order degenerate Bernoulli numbers. In the rest of this section, we recall the facts that are needed throughout this paper.\par

For any nonzero $\lambda\in\mathbb{R}$, the degenerate exponentials are defined by
\begin{equation}
e_{\lambda}^{x}(t)=\sum_{k=0}^{\infty}\frac{(x)_{k,\lambda}}{k!}t^{k}=(1+\lambda t)^{\frac{x}{\lambda}},\quad (\mathrm{see}\ [8-12]),\label{1}
\end{equation}
\begin{equation}
(x)_{0,\lambda}=1,\ (x)_{n,\lambda}=x(x-\lambda)(x-2\lambda)\cdots(x-(n-1)\lambda),\ (n\ge 1).\label{2}
\end{equation}
When $x=1$, $e_{\lambda}(t)=e_{\lambda}^{1}(t)=(1+\lambda t)^{\frac{1}{\lambda}}$, (see [7,13]).\par
In [3,4], Carlitz considered the degenerate Bernoulli polynomials of order $\alpha$ given by
\begin{equation}
\bigg(\frac{t}{e_{\lambda}(t)-1}\bigg)^{\alpha}e_{\lambda}^{x}(t)=\sum_{n=0}^{\infty}\beta_{n,\lambda}^{(\alpha)}(x)\frac{t^{n}}{n!}. \label{3}
\end{equation}
When $x=0$, $\beta_{n,\lambda}^{(\alpha)}=\beta_{n,\lambda}^{(\alpha)}(0)$ are called the degenerate Bernoulli numbers of order $\alpha$.  \par
\noindent In particular, for $\alpha=1$, $\beta_{n,\lambda}(x)=\beta_{n,\lambda}^{(1)}(x)$ are called the degenerate Bernoulli polynomials. \par
For $n\ge 0$, the degenerate Stirling numbers of the first kind are defined by
\begin{equation}
(x)_{n}=\sum_{k=0}^{n}S_{1,\lambda}(n,k)(x)_{k,\lambda},\quad (\mathrm{see}\ [8]),\label{4}
\end{equation}
where
\begin{displaymath}
(x)_{0}=1,\quad (x)_{n}=x(x-1)\cdots(x-n+1),\quad (n\ge 1),\quad (\mathrm{see}\ [1-19]).
\end{displaymath}
As the inversion formula of \eqref{4}, the degenerate Stirling numbers of the second kind are defined by
\begin{equation}
(x)_{n,\lambda}=\sum_{k=0}^{n}S_{2,\lambda}(n,k)(x)_{k},\quad (n\ge 0),\quad (\mathrm{see}\ [8]).\label{5}
\end{equation} \par
The degenerate logarithm $\log_{\lambda}(t)$ is the compositional inverse of $e_{\lambda}(t)$ satisfying $e_{\lambda}(\log_{\lambda}(t))=\log_{\lambda}(e_{\lambda}(t))$.
Then we have
\begin{equation}
\log_{\lambda}(1+t)=\frac{1}{\lambda}\big((1+t)^{\lambda}-1\big)=\sum_{n=1}^{\infty}\frac{\lambda^{n-1}(1)_{n,\frac{1}{\lambda}}}{n!}t^{n},\quad (\mathrm{see}\ [8]). \label{6}
\end{equation}
From \eqref{4} and \eqref{5}, we note that
\begin{equation}
\frac{1}{k!}\big(\log_{\lambda}(1+t)\big)^{k}=\sum_{n=k}^{\infty}S_{1,\lambda}(n,k)\frac{t^{n}}{n!},\quad (k\ge 0),\label{6-1}
\end{equation}
and
\begin{equation}
	\frac{1}{k!}\big(e_{\lambda}(t)-1\big)^{k}=\sum_{n=k}^{\infty}S_{2,\lambda}(n,k)\frac{t^{n}}{n!},\quad (\mathrm{see}\ [8]). \label{7}
\end{equation}
Let $\displaystyle f(t)=\sum_{k=0}^{\infty}a_{k}t^{k}\in\mathbb{C}[\![t]\!]\displaystyle$. For $n\ge 0$, the operator $[t^{n}]$ is defined by
\begin{equation}
[t^{n}]f(t)=a_{n},\quad (n\ge 0),\quad (\mathrm{see}\ [16]).\label{8}	
\end{equation}

\section{$r$-truncated degenerate Stirling numbers}
For $r\in\mathbb{N}$, we consider the $r$-truncated degenerate Stirling numbers of the second kind given by
\begin{equation}
\frac{1}{k!}\bigg(e_{\lambda}(t)-\sum_{l=0}^{r-1}\frac{(1)_{l,\lambda}}{l!}t^{l}\bigg)^{k}=\sum_{n=kr}^{\infty}S_{2,\lambda}^{[r]}(n,kr)\frac{t^{n}}{n!},\quad (k\ge 0). 	\label{9}
\end{equation}
We agree that $S_{2,\lambda}^{[r]}(n,kr)=0$, for $0 \le n <kr$. Note that $S_{2,\lambda}^{[1]}(n,k)=S_{2,\lambda}(n,k),\quad (n,k\ge 0)$. \par
From \eqref{9}, we note that
\begin{align}
&\frac{1}{k!}\bigg(e_{\lambda}(t)-\sum_{l=0}^{r-1}\frac{(1)_{l,\lambda}}{l!}t^{l}\bigg)^{k}=\frac{1}{k!}	\bigg(\sum_{l=r}^{\infty}\frac{(1)_{k,\lambda}}{l!}t^{l}\bigg)^{k}\label{10}\\
&=\frac{1}{k!}\sum_{n=kr}^{\infty}\bigg(\sum_{\substack{l_{1}+l_{2}+\cdots+l_{k}=n\\ l_{i}\ge r}}\frac{n!(1)_{l_{1},\lambda} (1)_{l_{2},\lambda}\cdots (1)_{l_{k},\lambda}}{l_{1}!l_{2}!\cdots l_{k}!}\bigg)\frac{t^{n}}{n!}. \nonumber
\end{align}
Thus, by \eqref{8} and \eqref{10}, we obtain the following theorem.
\begin{theorem}
For $n,k\ge 0$ with $n\ge kr$, we have
\begin{displaymath}
S_{2,\lambda}^{[r]}(n,kr)=\frac{1}{k!}\sum_{\substack{l_{1}+l_{2}+\cdots+l_{k}=n\\ l_{i}\ge r}}\frac{n!(1)_{l_{1},\lambda} (1)_{l_{2},\lambda}\cdots (1)_{l_{k},\lambda}}{l_{1}!l_{2}!\cdots l_{k}!}.
\end{displaymath}
\end{theorem}
From the binomial expansion, we note that
\begin{align}
&\frac{1}{k!}\bigg(e_{\lambda}(t)-\sum_{l=0}^{r-1}\frac{(1)_{l,\lambda}}{l!}t^{l}\bigg)^{k}\label{11}\\
&=\frac{1}{k!}\sum_{m=0}^{k}\binom{k}{m}e_{\lambda}^{k-m}(t)(-1)^{m}\bigg(\sum_{l=0}^{r-1}\frac{(1)_{l,\lambda}}{l!}t^{l}\bigg)^{m}\nonumber \\
&=\frac{1}{k!}\sum_{m=0}^{k}\binom{k}{m}\sum_{j=0}^{\infty}(k-m)_{j,\lambda}\frac{t^{j}}{j!}(-1)^{m}\sum_{l_{1},l_{2},\dots,l_{m}=0}^{r-1}\frac{(1)_{l_{1},\lambda}(1)_{l_{2},\lambda}\cdots (1)_{l_{m},\lambda}}{l_{1}!l_{2}!\cdots l_{m}!}t^{l_{1}+\cdots+l_{m}}\nonumber \\
&=\sum_{n=0}^{\infty}\bigg(\frac{1}{k!}\sum_{m=0}^{k}\binom{k}{m}(-1)^{m}\sum_{l_{1},l_{2},\dots,l_{m}=0}^{r-1} \frac{n!(1)_{l_{1},\lambda} (1)_{l_{2},\lambda}\cdots (1)_{l_{m},\lambda}(k-m)_{n-l_{1}-l_{2}-\cdots-l_{m,\lambda}}}{l_{1}!l_{2}!\cdots l_{m}!(n-l_{1}-\cdots-l_{m})!}
\bigg)\frac{t^{n}}{n!}.\nonumber
\end{align}
Therefore, by comparing the coefficients on both sides of \eqref{11}, we obtain the following theorem.
\begin{theorem}
For $n,k\ge 0$, we have
\begin{align*}
& \frac{1}{k!}\sum_{m=0}^{k}\binom{k}{m}(-1)^{m}\sum_{l_{1},l_{2},\dots,l_{m}=0}^{r-1} \frac{n!(1)_{l_{1},\lambda} (1)_{l_{2},\lambda}\cdots (1)_{l_{m},\lambda}(k-m)_{n-l_{1}-l_{2}-\cdots-l_{m,\lambda}}}{l_{1}!l_{2}!\cdots l_{m}!(n-l_{1}-\cdots-l_{m})!}\\
&=\left\{\begin{array}{ccc} S_{2,\lambda}^{[r]}(n,kr), & \textrm{ if $n\ge kr$,} \\
0, & \textrm{if $0 \le n<kr$.}
\end{array}\right.
\end{align*}
\end{theorem}
When $k=1$ in \eqref{9}, we have
\begin{align}
&\sum_{l=r}^{\infty}\frac{(1)_{l,\lambda}}{l!}t^{l}=\bigg(e_{\lambda}(t)-\sum_{l=1}^{r-1}\frac{(1)_{l,\lambda}}{l!}t^{l}\bigg)	\label{12}\\
&\quad =\sum_{n=r}^{\infty}S_{2,\lambda}^{[r]}(n,r)\frac{t^{n}}{n!}=\sum_{n=0}^{\infty}S_{2,\lambda}^{[r]}(n+r,r)\frac{t^{n+r}}{(n+r)!} \nonumber \\
&=t^{r}\sum_{n=0}^{\infty}S_{2,\lambda}^{[r]}(n+r,r)\frac{n!}{(n+r)!}\frac{t^{n}}{n!}.\nonumber
\end{align}
For $k=2$, we have
\begin{align}
&\bigg(\sum_{l=r}^{\infty}\frac{(1)_{l,\lambda}}{l!}t^{l}\bigg)^{2} \label{13} \\
&=t^{2r}\bigg(\sum_{j=0}^{\infty}S_{2,\lambda}^{[r]}(j+r,r)\frac{t^{j}}{(j+r)!}\bigg)\bigg(\sum_{l=0}^{\infty}S_{2,\lambda}^{[r]}(l+r,r)\frac{t^{l}}{(l+r)!}\bigg) \nonumber \\
&=t^{2r}\sum_{n=0}^{\infty}\bigg(\sum_{j=0}^{n}S_{2,\lambda}^{[r]}(j+r,r)S_{2,\lambda}^{[r]}(n-j+r,r)\frac{n!}{(n-j+r)!(j+r)!}\bigg)\frac{t^{n}}{n!}. \nonumber
\end{align}
Continuing this process, we have
\begin{align}
&\bigg(\sum_{l=r}^{\infty}\frac{(1)_{l,\lambda}}{l!}t^{l}\bigg)^{k} \label{14}	\\
&=t^{kr}\sum_{n=0}^{\infty}\bigg(\sum_{j_{1}+j_{2}+\cdots+j_{k}=n}\frac{n!S_{2,\lambda}^{[r]}(j_{1}+r,r) S_{2,\lambda}^{[r]}(j_{2}+r,r)\cdots S_{2,\lambda}^{[r]}(j_{k}+r,r) }{(j_{1}+r)!(j_{2}+r)!\cdots(j_{k}+r)!} \bigg)\frac{t^{n}}{n!}.\nonumber
\end{align}
On the other hand, by \eqref{9}, we get
\begin{align}
&\bigg(\sum_{l=r}^{\infty}\frac{t^{l}}{l!}(1)_{l,\lambda}\bigg)^{k}=\bigg(e_{\lambda}(t)-\sum_{l=1}^{r-1}\frac{(1)_{l,\lambda}}{l!}t^{l}\bigg)^{k}\label{15} \\
&=k!\sum_{n=kr}^{\infty}S_{2,\lambda}^{[r]}(n,kr)\frac{t^{n}}{n!}=k!\sum_{n=0}^{\infty}S_{2,\lambda}^{[r]}(n+kr,kr)\frac{t^{n+kr}}{(n+kr)!} \nonumber \\
&=t^{kr}\sum_{n=0}^{\infty}\frac{k!S_{2,\lambda}^{[r]}(n+kr,kr)}{(n+kr)!}t^{n}.\nonumber
\end{align}
Therefore, by \eqref{14} and \eqref{15}, we obtain the following theorem.
\begin{theorem}
For $n,k\ge 0$, we have
\begin{align*}
	&\frac{k!}{(n+kr)!}S_{2,\lambda}^{[r]}(n+kr,kr)\\
	&=\sum_{j_{1}+j_{2}+\cdots+j_{k}=n}\frac{S_{2,\lambda}^{[r]}(j_{1}+r,r) S_{2,\lambda}^{[r]}(j_{2}+r,r)\cdots S_{2,\lambda}^{[r]}(j_{k}+r,r)}{(j_{1}+r)!(j_{2}+r)!\cdots(j_{k}+r)!}.
\end{align*}	
\end{theorem}
Let us consider the $r$-truncated degenerate Bernoulli polynomials of order $\alpha$ given by
\begin{equation}
\frac{t^{\alpha r}}{\big(e_{\lambda}(t)-\sum_{l=0}^{r-1}\frac{t^{l}}{l!}(1)_{l,\lambda}\big)^{\alpha}}e_{\lambda}^{x}(t)	=\sum_{n=0}^{\infty}\beta_{n,\lambda}^{[r-1,\alpha]}(x)\frac{t^{n}}{n!}. \label{16}
\end{equation}
When $x=0$, $\beta_{n,\lambda}^{[r-1,\alpha]}=\beta_{n,\lambda}^{[r-1,\alpha]}(0)$ are called the $r$-truncated degenerate Bernoulli numbers of order $\alpha$. \par
Note that
\begin{align}
t^{\alpha r}&=\bigg(e_{\lambda}(t)-\sum_{l=0}^{r-1}\frac{t^{l}}{l!}(1)_{l,\lambda}\bigg)^{\alpha}\sum_{l=0}^{\infty}\beta_{l,\lambda}^{[r-1,\alpha]}\frac{t^{l}}{l!} \nonumber  \\
&=\alpha!\sum_{j=\alpha r}^{\infty} S_{2,\lambda}^{[r]}(j,\alpha r)\frac{t^{j}}{j!}\sum_{l=0}^{\infty}\beta_{l,\lambda}^{[r-1,\alpha]}\frac{t^{l}}{l!} \label{17} \\
&=\alpha!\sum_{n=\alpha r}^{\infty}\bigg(\sum_{l=0}^{n-\alpha r}\frac{S_{2,\lambda}^{[r]}(n-l,\alpha r)n!}{(n-l)!l!}\beta_{l,\lambda}^{[r-1,\alpha]}\bigg)\frac{t^{n}}{n!} \nonumber \\
&=\sum_{n=\alpha r}^{\infty}\bigg(\alpha!\sum_{l=0}^{n-\alpha r}\binom{n}{l}S_{2,\lambda}^{[r]}(n-l,\alpha r)\beta_{l,\lambda}^{[r-1,\alpha]}\bigg)\frac{t^{n}}{n!},\nonumber
\end{align}
where $\alpha$ is a positive integer. \par
\noindent Thus, by \eqref{17}, we get
\begin{displaymath}
	\sum_{l=0}^{n-\alpha r}\binom{n}{l}S_{2,\lambda}^{[r]}(n-l,\alpha r)\beta_{l,\lambda}^{[r-1,\alpha]}=\left\{\begin{array}{ccc}
	\frac{(\alpha r)!}{\alpha!}, & \textrm{if $n=\alpha r$,} \\
	0, & \textrm{if $ n > \alpha r$.}
\end{array}\right.
\end{displaymath}

It is well known that the partial Bell polynomials are defined by
\begin{displaymath}
	\frac{1}{k!}\bigg(\sum_{l=1}^{\infty}x_{l}\frac{t^{l}}{l!}\bigg)^{k}=\sum_{n=k}^{\infty}B_{n,k}(x_{1},x_{2},\dots,x_{n-k+1})\frac{t^{n}}{n!},
\end{displaymath}
where $k$ is a nonnegative integer. \par
\noindent Thus, we note that
\begin{align*}
&B_{n,k}(x_{1},x_{2},\dots,x_{n-k+1})\\
&=\sum_{\substack{l_{1}+l_{2}+\cdots+l_{n-k+1}=k\\ l_{1}+2l_{2}+\cdots+(n-k+1)l_{n-k+1}=n}}\frac{n!}{l_{1}!l_{2}!\cdots l_{n-k+1}!}\bigg(\frac{x_{1}}{1!}\bigg)^{l_{1}}\bigg(\frac{x_{2}}{2!}\bigg)^{l_{2}}\cdots \bigg(\frac{x_{n-k+1}}{(n-k+1)!}\bigg)^{l_{n-k+1}}\\
&=n!\sum_{\Lambda_{n}^{k}}\prod_{j=1}^{n-k+1}\frac{1}{k_{j}!}\bigg(\frac{x_{j}}{j!}\bigg)^{k_{j}},
\end{align*}
where
\begin{displaymath}
	\Lambda_{n}^{k}=\{(k_{1},k_{2},\dots,k_{n-k+1})\ |\ k_{1}+k_{2}+\cdots+k_{n-k+1}=k,\ k_{1}+2k_{2}+\cdots+(n-k+1)k_{n-k+1}=n\}.
\end{displaymath} \par
We observe that
\begin{align}
&\frac{1}{1+\frac{x_{1}}{1!}(1)_{1,\lambda}t+\frac{x_{2}}{2!}(1)_{2,\lambda}t^{2}+\cdots+\frac{x_{k}}{k!}(1)_{k,\lambda}t^{k}+\cdots } \label{18} \\
&=\sum_{k=0}	^{\infty}(-1)^{k}k!\frac{1}{k!}\bigg(\sum_{l=1}^{\infty}\frac{x_{l}(1)_{l,\lambda}}{l!}t^{l}\bigg)^{k} \nonumber \\
&=1+\sum_{k=1}^{\infty}(-1)^{k}k!\sum_{n=k}^{\infty}B_{n,k}(x_{1}(1)_{1,\lambda},x_{2}(1)_{2,\lambda},\dots,x_{n-k+1}(1)_{n-k+1})\frac{t^{n}}{n!} \nonumber \\
&=1+\sum_{n=1}^{\infty}\bigg(\sum_{k=1}^{n}(-1)^{k}k!B_{n,k}(x_{1}(1)_{1,\lambda},x_{2}(1)_{2,\lambda},\dots,(x)_{n-k+1,\lambda})\bigg)\frac{t^{n}}{n!}. \nonumber
\end{align}
We denote \eqref{18} by
\begin{equation}
\frac{1}{\sum_{n=0}^{\infty}\frac{(1)_{n,\lambda}}{n!}x_{n}t^{n}}=\sum_{n=0}^{\infty}K_{n,\lambda}(x_{1},x_{2},\dots)\frac{t^{n}}{n!}. \label{19}	
\end{equation}
Then, by \eqref{18} and \eqref{19}, we get
\begin{align*}
K_{n,\lambda}(x_{1},x_{2},\dots)&=\sum_{k=1}^{n}(-1)^{k}	k!B_{n,k}(x_{1}(1)_{1,\lambda},x_{2}(1)_{2,\lambda},\dots,x_{n-k+1}(1)_{n-k+1,\lambda}\big), \quad(n \ge 1),\\
&K_{0,\lambda}(x_{1},x_{2},\dots)=1,
\end{align*} \par
\noindent Note that, from \eqref{18} and \eqref{19}, we have
\begin{align*}
	K_{n,\lambda}(1,1,\dots)&=\sum_{k=0}^{n}(-1)^{k}k!B_{n,k}\big((1)_{1,\lambda},(1)_{2,\lambda},\dots,(1)_{n-k+1,\lambda}\big) \\
	&=\sum_{k=0}^{n}(-1)^{k}k!S_{2,\lambda}(n,k).
\end{align*} \par
Taking $\alpha=1$ in \eqref{16}, we have
\begin{equation}
	\frac{t^{r}}{e_{\lambda}(t)-\sum_{l=0}^{r-1}\frac{(1)_{l,\lambda}}{l!}t^{l}}e_{\lambda}^{x}(t)=\sum_{n=0}^{\infty}\beta_{n,\lambda}^{[r-1,1]}(x)\frac{t^{n}}{n!}. \label{20}
\end{equation}
When $r=1$, $\beta_{n,\lambda}^{[0,1]}(x)=\beta_{n,\lambda}(x)$ are the degenerate Bernoulli polynomials. \par
\noindent From \eqref{20}, we note that
\begin{align}
e_{\lambda}^{x}(t)&=\frac{1}{t^{r}}\bigg(e_{\lambda}(t)-\sum_{l=1}^{r-1}\frac{(1)_{l,\lambda}}{l!}t^{l}\bigg)\sum_{l=0}^{\infty}\beta_{l,\lambda}^{[r-1,1]}(x)\frac{t^{l}}{l!} \label{21}\\
&=\frac{1}{t^{r}}\sum_{j=r}^{\infty}\frac{(1)_{j,\lambda}}{j!}t^{j}\sum_{l=0}^{\infty}\beta_{l,\lambda}^{[r-1,1]}(x)\frac{t^{l}}{l!} \nonumber \\
&=\sum_{j=0}^{\infty}\frac{(1)_{j+r,\lambda}}{(j+r)!}t^{j}\sum_{l=0}^{\infty}\beta_{l,\lambda}^{[r-1,1]}(x)\frac{t^{l}}{l!} \nonumber \\
&=\sum_{n=0}^{\infty}\bigg(\sum_{j=0}^{n}\binom{n}{j}\frac{j!(1)_{j+r,\lambda}}{(j+r)!}\beta_{n-j,\lambda}^{[r-1,1]}(x)\bigg)\frac{t^{n}}{n!}.\nonumber
\end{align}
By \eqref{1} and \eqref{21}, we get
\begin{equation}
(x)_{n,\lambda}=\sum_{j=0}^{n}\binom{n}{j}\frac{j!}{(j+r)!}(1)_{j+r,\lambda}\beta_{n-j,\lambda}^{[r-1,1]}(x),\quad (n\ge 0), \label{22}	
\end{equation}
In particular, for $n=0$, $\beta_{0,\lambda}^{[r-1,1]}=\frac{r!}{(1)_{r,\lambda}}$. \par
\noindent By \eqref{22}, we get
\begin{align*}
\beta_{1,\lambda}^{[r-1,1]}&=\frac{r!}{(1)_{r,\lambda}}\bigg(x-\frac{1-r \lambda}{1+r}\bigg),\\
\beta_{2,\lambda}^{[r-1,1]}&=\frac{r!}{(1)_{r,\lambda}}\bigg\{x(x-\lambda)+\frac{2x(1-r\lambda)}{1+r}-\frac{2(1-r\lambda)^{2}}{(1+r)^{2}} \\
&\qquad\qquad+\frac{2}{(1+r)(2+r)}(1-r\lambda)(1-(r+1)\lambda)\bigg\}.
\end{align*} \par
From \eqref{6}, we note that
\begin{equation}
\frac{1}{t}\log_{\lambda}(1+t)=\frac{1}{t}\sum_{n=1}^{\infty}\frac{\lambda^{n-1}(1)_{n,\frac{1}{\lambda}}}{n!}t^{n}	=\sum_{n=0}^{\infty}\frac{\lambda^{n}(1)_{n+1,\frac{1}{\lambda}}}{(n+1)!}t^{n}.\label{24}
\end{equation}
Replacing $t$ by $e_{\lambda}(t)-1$ in \eqref{24}, we get
\begin{align}
\sum_{n=0}^{\infty}\beta_{n,\lambda}\frac{t^{n}}{n!}&=\frac{t}{e_{\lambda}(t)-1}=\sum_{k=0}^{\infty}\frac{\lambda^{k}(1)_{k+1,\frac{1}{\lambda}}}{(k+1)}\frac{1}{k!}\big(e_{\lambda}(t)-1\big)^{k} \label{25} \\
&=\sum_{n=0}^{\infty}\bigg(\sum_{k=0}^{n}\frac{\lambda^{k}(1)_{k+1,\frac{1}{\lambda}}}{k+1}S_{2,\lambda}(n,k)\bigg)\frac{t^{n}}{n!}.\nonumber
\end{align}
Comparing the coefficients on both sides of \eqref{25}, we have
\begin{equation}
\beta_{n,\lambda}=\sum_{k=0}^{n}\frac{\lambda^{k}(1)_{k+1,\frac{1}{\lambda}}}{k+1}S_{2,\lambda}(n,k),\quad (n\ge 0). \label{26}	
\end{equation} \par
Replacing $t$ by $\log_{\lambda}(1+t)$ in \eqref{3}, with $\alpha=1, x=0$, we have
\begin{align}
\sum_{n=0}^{\infty}\frac{\lambda^{n}(1)_{n+1,\frac{1}{\lambda}}}{(n+1)!}t^{n}&=\frac{1}{t}\log_{\lambda}(1+t)=\sum_{k=0}^{\infty}\beta_{k,\lambda}\frac{1}{k!}\big(\log_{\lambda}(1+t)\big)^{k}\label{27}\\
&=\sum_{k=0}^{\infty}\beta_{k,\lambda}\sum_{n=k}^{\infty}S_{1,\lambda}(n,k)\frac{t^{n}}{n!} \nonumber \\
&=\sum_{n=0}^{\infty}\bigg(\sum_{k=0}^{n}\beta_{k,\lambda}S_{1,\lambda}(n,k)\bigg)\frac{t^{n}}{n!}.\nonumber
\end{align}
By comparing the coefficients on both sides of \eqref{27}, we get
\begin{equation}
\sum_{k=0}^{n}\beta_{k,\lambda}S_{1,\lambda}(n,k)=\frac{\lambda^{n}(1)_{n+1,\frac{1}{\lambda}}}{n+1},\quad (n\ge 0). \label{28}	
\end{equation}
Therefore, by \eqref{26} and \eqref{28}, we obtain the following theorem.
\begin{theorem}
For $n\ge 0$, we have
\begin{displaymath}
\beta_{n,\lambda}=\sum_{k=0}^{n}\frac{\lambda^{k}(1)_{k+1,\frac{1}{\lambda}}}{k+1}S_{2,\lambda}(n,k),
\end{displaymath}
and
\begin{displaymath}
\frac{(1)_{n+1,\frac{1}{\lambda}}}{n+1}\lambda^{n}=\sum_{k=0}^{n}\beta_{k,\lambda}S_{1,\lambda}(n,k).
\end{displaymath}
\end{theorem} \par
From Theorem 4, we note that
\begin{align}
&\sum_{j=0}^{n}S_{2,\lambda}^{[2]}(n-j+k,2k)\binom{n+k}{j}\beta_{j,\lambda}\label{29}\\
&=(n+k)!\sum_{j=0}^{n}\frac{S_{2,\lambda}^{[2]}(n-j+k,2k)}{(n-j+k)!}\frac{\beta_{j,\lambda}}{j!}\nonumber\\
&=(n+k)!\sum_{j=0}^{n}[t^{n-j+k}]\frac{(e_{\lambda}(t)-1-t)^{k}}{k!}[t^{j}]\frac{t}{e_{\lambda}(t)-1} \nonumber \\
&=(n+k)!\sum_{j=0}^{n}[t^{n-j}]\frac{(e_{\lambda}(t)-1-t)^{k}}{k!t^{k}}[t^{j}]\frac{t}{e_{\lambda}(t)-1} \nonumber \\
&=(n+k)![t^{n}]\frac{(e_{\lambda}(t)-1-t)^{k}t}{k!t^{k}(e_{\lambda}(t)-1)}\nonumber\\
 &=(n+k)![t^{n}]\sum_{j=0}^{k}\binom{k}{j}\frac{(-1)^{k-j}(e_{\lambda}(t)-1)^{j}t^{k-j}t}{k!t^{k}(e_{\lambda}(t)-1)}\nonumber \\
 &=\frac{(-1)^{k}(n+k)!}{k!}[t^{n}]\frac{t}{e_{\lambda}(t)-1}+(n+k)!\sum_{j=1}^{k}(-1)^{k-j}[t^{n}]\frac{(e_{\lambda}(t)-1)^{j-1}}{j!(k-j)!t^{j-1}}\nonumber \\
 &=\frac{(-1)^{k}(n+k)!}{k!}\frac{\beta_{n,\lambda}}{n!}+(n+k)!\sum_{j=1}^{k}\frac{(-1)^{k-j}}{j!(k-j)!}[t^{n}]\frac{(e_{\lambda}(t)-1)^{j-1}}{t^{j-1}}\nonumber \\
 &=(-1)^{k}\binom{n+k}{k}\beta_{n,\lambda}+(n+k)!\sum_{j=1}^{k}\frac{(-1)^{k-j}(j-1)!}{j!(k-j)!}[t^{n}] \sum_{l=0}^{\infty}\frac{S_{2,\lambda}(l+j-1,j-1)}{(l+j-1)!}t^{l}\nonumber \\
 &=(-1)^{k}\binom{n+k}{k}\beta_{n,\lambda}+(n+k)!\sum_{j=1}^{k}\frac{(-1)^{k-j}}{j(k-j)!}\frac{S_{2,\lambda}(n+j-1,j-1)}{(n+j-1)!}.\nonumber
\end{align}
Therefore, by \eqref{29}, we obtain the following theorem.
\begin{theorem}
For $n, k \ge 0$, we have
\begin{align*}
&\sum_{j=0}^{n}S_{2,\lambda}^{[2]}(n-j+k,2k)\binom{n+j}{k}\beta_{j,\lambda}\\
&=(-1)^{n}\binom{n+k}{k}\beta_{n,\lambda}+(n+k)!\sum_{j=1}^{k}\frac{(-1)^{k-j}}{j(k-j)!(n+j-1)!}S_{2,\lambda}(n+j-1,j-1).
\end{align*}
\end{theorem}
Observe that, for any formal power series $f(t)=\sum_{i=0}^{\infty}a_{i}t^{i}$, we have
\begin{equation}
[t^{n}]f^{\prime}(t)=(n+1)a_{n+1}=(n+1)[t^{n+1}]f(t).\label{30}	
\end{equation}
By making use of \eqref{30}, we have
\begin{align}
&\sum_{j=0}^{n}\binom{n+k-1}{j}S_{2,\lambda}^{[2]}(n-j+k,2k)\beta_{j,\lambda} \label{31}\\
&=(n+k-1)!\sum_{j=0}^{n}\frac{S_{2,\lambda}^{[2]}(n-j+k,2k)}{(n-j+k-1)!}\frac{\beta_{j,\lambda}}{j!} \nonumber \\
&=(n+k-1)!\sum_{j=0}^{n}(n-j+k)[t^{n-j+k}]\frac{(e_{\lambda}(t)-1-t)^{k}}{k!}[t^{j}]\frac{t}{e_{\lambda}(t)-1} \nonumber \\
&=(n+k-1)!\sum_{j=0}^{n}[t^{n-j+k-1}]\frac{d}{dt}\bigg(\frac{(e_{\lambda}(t)-1-t)^{k}}{k!}\bigg)[t^{j}]\frac{t}{e_{\lambda}(t)-1}\nonumber \\
&=(n+k-1)!\sum_{j=0}^{n}[t^{n-j+k-1}]\frac{(e_{\lambda}(t)-1-t)^{k-1}(e_{\lambda}^{1-\lambda}(t)-1)}{(k-1)!}[t^{j}]\frac{t}{e_{\lambda}(t)-1} \nonumber \\
&=(n+k-1)!\sum_{j=0}^{n}[t^{n-j+k-1}]\frac{(e_{\lambda}(t)-1-t)^{k-1}(e_{\lambda}(t)-1-\lambda t)}{(k-1)!(1+\lambda t)}[t^{j}]\frac{t}{e_{\lambda}(t)-1} \nonumber \\
&=(n+k-1)!\sum_{j=0}^{n}[t^{n-j}]\frac{(e_{\lambda}(t)-1-t)^{k-1}(e_{\lambda}(t)-1-\lambda t)}{t^{k-1}(k-1)!(1+\lambda t)}[t^{j}]\frac{t}{e_{\lambda}(t)-1} \nonumber \\
&=(n+k-1)![t^{n}]\frac{(e_{\lambda}(t)-1-t)^{k-1}}{t^{k-2}(k-1)!(1+\lambda t)} \nonumber \\
&\qquad -(n+k-1)!\lambda [t^{n}]\frac{(e_{\lambda}(t)-1-t)^{k-1}}{t^{k-2}(k-1)!(1+\lambda t)}\frac{t}{e_{\lambda}(t)-1}\nonumber \\
&=(n+k-1)![t^{n}]\bigg\{\sum_{m=k}^{\infty}\sum_{l=k}^{m}\frac{S_{2,\lambda}^{[2]}(l+k-2,2k-2)}{(l+k-2)!}(-\lambda)^{m-l}t^{m}\nonumber \\
&\qquad -\lambda\sum_{m=k}^{\infty}\sum_{j=k}^{m}\sum_{l=k}^{j}\frac{S_{2,\lambda}^{[2]}(l+k-2,2k-2)}{(l+k-2)!}(-\lambda)^{j-l}\frac{\beta_{m-j,\lambda}}{(m-j)!}t^{m}\bigg\} \nonumber \\
&=(n+k-1)!\sum_{l=k}^{n}\frac{S_{2,\lambda}^{[2]}(l+k-2,2k-2)}{(l+k-2)!}(-\lambda)^{n-l}\nonumber\\
&\qquad-\lambda (n+k-1)!\sum_{j=k}^{n}\sum_{l=k}^{j}\frac{S_{2,\lambda}^{[2]}(l+k-2,2k-2)\beta_{n-j,\lambda}(-\lambda)^{j-l}}{(l+k-2)!(n-j)!}\nonumber.
\end{align}
Therefore, by \eqref{31}, we obtain the following theorem.
\begin{theorem}
For $n,k\ge 1$, with $n \ge k$, we have
\begin{align*}
& \sum_{j=0}^{n}\binom{n+k-1}{j}S_{2,\lambda}^{[2]}(n-j+k,2k)\beta_{j,\lambda}\\
&=(n+k-1)!\bigg\{\sum_{l=k}^{n}\frac{S_{2,\lambda}^{[2]}(l+k-2,2k-2)}{(l+k-2)!}(-\lambda)^{n-l}\nonumber\\
&\qquad\qquad\qquad +\sum_{j=k}^{n}\sum_{l=k}^{j}\frac{S_{2,\lambda}^{[2]}(l+k-2,2k-2)\beta_{n-j,\lambda}(-\lambda)^{j-l+1}}{(l+k-2)!(n-j)!}\bigg\}.
\end{align*}
\end{theorem}
Now, we observe that
\begin{align}
&\sum_{j=0}^{n}S_{2,\lambda}^{[2]}(n-j+k,2k)\binom{n+k}{j}\beta_{j,\lambda}^{(k)} \label{32}\\
&=(n+k)!\sum_{j=0}^{n}\frac{S_{2,\lambda}^{[2]}(n-j+k,2k)}{(n-j+k)!}\frac{\beta_{j,\lambda}^{(k)}}{j!}\nonumber \\
&=(n+k)!\sum_{j=0}^{n}[t^{n-j+k}]\frac{(e_{\lambda}(t)-1-t)^{k}}{k!} [t^{j}]\bigg(\frac{t}{e_{\lambda}(t)-1}\bigg)^{k} \nonumber \\
&=\frac{(n+k)!}{k!}\sum_{j=0}^{n}[t^{n-j}]\frac{(e_{\lambda}(t)-1-t)^{k}}{t^{k}}[t^{j}] \bigg(\frac{t}{e_{\lambda}(t)-1}\bigg)^{k} \nonumber \\
&=\frac{(n+k)!}{k!}[t^{n}]\bigg(1-\frac{t}{e_{\lambda}(t)-1}\bigg)^{k}\nonumber \\
&=\frac{(n+k)!}{k!}\sum_{j=0}^{k}\binom{k}{j}(-1)^{k-j}[t^{n}]\bigg(\frac{t}{e_{\lambda}(t)-1}\bigg)^{k-j} \nonumber \\
&=\frac{(n+k)!}{k!}\sum_{j=0}^{k}\binom{k}{j}(-1)^{k-j}\frac{\beta_{n,\lambda}^{(k-j)}}{n!}. \nonumber
\end{align}
Therefore, by \eqref{32}, we obtain the following theorem.
\begin{theorem}
For $n,k\ge 0$, we have
\begin{displaymath}
\sum_{j=0}^{n}S_{2,\lambda}^{[2]}(n-j+k,2k)\binom{n+k}{j}\beta_{j,\lambda}^{(k)}=\binom{n+k}{k}\sum_{j=0}^{k}\binom{k}{j}(-1)^{k-j}\beta_{n,\lambda}^{(k-j)}.
\end{displaymath}
\end{theorem}
From Theorems 4 and 5, we have
\begin{align}
&\sum_{j=0}^{n}\binom{n+k}{j+k}S_{2,\lambda}^{[3]}(j+k,3k)\beta_{n-j,\lambda}=(n+k)!\sum_{j=0}^{n}\frac{S_{2,\lambda}^{[3]}(j+k,3k)}{(j+k)!} \frac{\beta_{n-j,\lambda}}{(n-j)!}\label{33}\\
&=(n+k)!\sum_{j=0}^{n}[t^{j+k}]\frac{(e_{\lambda}(t)-1-t-(1)_{2,\lambda}t^{2}/2)^{k}}{k!}[t^{n-j}]\frac{t}{e_{\lambda}(t)-1} \nonumber \\
&=(n+k)!\sum_{j=0}^{n}[t^{j}]	\frac{(e_{\lambda}(t)-1-t-(1)_{2,\lambda}t^{2}/2)^{k}}{t^{k}k!}[t^{n-j}]\bigg(\frac{t}{e_{\lambda}(t)-1}\bigg)\nonumber \\
&=(n+k)![t^{n}]\frac{(e_{\lambda}(t)-1-t-(1)_{2,\lambda}t^{2}/2)^{k}t}{t^{k}k!(e_{\lambda}(t)-1)}\nonumber \\
&=(n+k)![t^{n}]\sum_{j=0}^{k}(-1)^{k-j}\binom{k}{j}\frac{(e_{\lambda}(t)-1)^{j}(t+(1-\lambda)\frac{t^{2}}{2})^{k-j}t}{k!(e_{\lambda}(t)-1)t^{k}}\nonumber
\end{align}
\begin{align*}
&=(n+k)![t^{n}]\frac{(-1)^{k}\big(1+(1-\lambda)\frac{t}{2}\big)^{k}t}{k!(e_{\lambda}(t)-1)} \nonumber \\
&\qquad +\frac{(n+k)!}{k!}[t^{n}]\sum_{j=1}^{k}(-1)^{k-j}\binom{k}{j}(e_{\lambda}(t)-1)^{j-1}t^{1-j}\bigg(1+(1-\lambda)\frac{t}{2}\bigg)^{k-j}\nonumber \\
&=(n+k)![t^{n}]\sum_{m=0}^{\infty}\bigg((-1)^{k}\sum_{j=0}^{m}\binom{k}{j}(1-\lambda)^{j}\frac{\beta_{m-j,\lambda}}{2^{j}(m-j)!k!}\bigg)t^{m}\nonumber \\
&\qquad +\frac{(n+k)!}{k!}[t^{n}]\sum_{j=1}^{k}(-1)^{k-j}\binom{k}{j}(e_{\lambda}(t)-1)^{j-1}\sum_{l=0}^{k-j}\binom{k-j}{l}(1-\lambda)^{l}\frac{t^{l-j+1}}{2^{l}}\nonumber \\
&=(-1)^{k}(n+k)!\sum_{j=0}^{n}\binom{k}{j}(1-\lambda)^{j}\frac{1}{2^{j}(n-j)!k!}\beta_{n-j,\lambda}\nonumber \\
&+(n+k)![t^{n}]\sum_{j=1}^{k}\frac{(-1)^{k-j}}{j(k-j)!}\sum_{l_{1}=j-1}^{\infty}S_{2,\lambda}(l_{1},j-1)\frac{t^{l_{1}}}{l_{1}!}\sum_{l=1}^{k-j}\binom{k-j}{l}\frac{(1-\lambda)^{l}}{2^{l}}t^{l-j+1}\nonumber \\
&=(-1)^{k}(n+k)!\sum_{j=0}^{n}\binom{k}{j}(1-\lambda)^{j}\frac{1}{2^{j}(n-j)!k!}\beta_{n-j,\lambda} \nonumber \\
&\qquad +(n+k)!\sum_{j=1}^{k}\frac{(-1)^{k-j}}{j(k-j)!}\sum_{l=1}^{n}\binom{k-j}{l}\frac{S_{2,\lambda}(n-l+j-1,j-1)(1-\lambda)^{l}}{2^{l}(n-l+j-1)!}.\nonumber
\end{align*}
Therefore, by \eqref{33}, we obtain the following theorem.
\begin{theorem}
For $n,k\ge 0$, we have
\begin{align}
&\sum_{j=0}^{n}\binom{n+k}{j+k}S_{2,\lambda}^{[3]}(j+k,3k)\beta_{n-j,\lambda} \\
 &=(-1)^{k}(n+k)!\sum_{j=0}^{n}\binom{k}{j}(1-\lambda)^{j}\frac{1}{2^{j}(n-j)!k!}\beta_{n-j,\lambda} \nonumber \\
&\qquad +(n+k)!\sum_{j=1}^{k}\frac{(-1)^{k-j}}{j(k-j)!}\sum_{l=1}^{n}\binom{k-j}{l}\frac{S_{2,\lambda}(n-l+j-1,j-1)(1-\lambda)^{l}}{2^{l}(n-l+j-1)!}\nonumber.
\end{align}
\end{theorem}

\begin{remark}
As the counterpart of \eqref{9}, we may consider the $r$-truncated degenerate Stirling numbers of the first kind given by
\begin{equation}
\frac{1}{k!}\bigg(\log_{\lambda}(1+t)-\sum_{l=1}^{r-1}\frac{(1)_{l,\frac{1}{\lambda}}\lambda^{l-1}}{l!}t^{l}\bigg)^{k}=\sum_{n=kr}^{\infty}S_{1,\lambda}^{[r]}(n,kr)\frac{t^{n}}{n!},\label{23}
\end{equation}
where $r$ is a positive integer. These numbers will be investigated in a forthcoming paper.
\end{remark}

\section{Conclusion}
In recent years, we have witnessed that many degenerate versions of quite a few special numbers and polynomials were investigated and some nice results were obtained by adopting various tools. \par
In this paper, we considered the $r$-truncated degenerate Stirling numbers of the second, which reduce to the degenerate Stirling numbers of the second for $r=1$, and studied by using generating functions their explicit expressions, some properties and related identities on those numbers, in connection with several other degenerate special numbers and polynomials. \par
As one of our future projects, we would like to continue this line of research, namely, to explore various degenerate versions of some special numbers and polynomials, and to find their applications to physics, science and engineering.

\vspace{0.1in}

\noindent{\bf{Funding}} \\
This work was supported by the Basic Science Research Program, the National
               Research Foundation of Korea,
                (NRF-2021R1F1A1050151).

\end{document}